\documentclass{amsart}
\usepackage{amssymb,amscd,verbatim}
\numberwithin{equation}{section}

\newtheorem{prop}{Proposition}
\newtheorem{thm}[prop]{Theorem}
\newtheorem{cor}[prop]{Corollary}
\newtheorem{lem}[prop]{Lemma}

\theoremstyle{definition}
\newtheorem{ex}[prop]{Example}
\newtheorem{rem}[prop]{Remark}

\begin{document}

\title[]{Affine type A Crystal Structure on Tensor Products of Rectangles,
Demazure characters, and Nilpotent Varieties}

\author{Mark Shimozono}
\address{Dept.\ of Mathematics \\ Virginia Tech \\ Blacksburg, VA}
\email{mshimo@math.vt.edu}

\begin{abstract}
Answering a question of Kuniba, Misra, Okado, Takagi, and Uchi-yama,
it is shown that certain Demazure characters of affine type A,
coincide with the graded characters of coordinate rings of closures
of conjugacy classes of nilpotent matrices.
\end{abstract}

\maketitle

\newcommand{\id}{\mathrm{id}}
\newcommand{\content}{\mathrm{content}}
\newcommand{\ch}{\mathrm{ch}}
\newcommand{\word}{\mathrm{word}}
\newcommand{\charge}{\mathrm{charge}}
\newcommand{\key}{\mathrm{Key}}
\newcommand{\Hom}{\mathrm{Hom}}
\newcommand{\PCL}{P_{cl}}
\newcommand{\PCLB}{\overline{P}_{cl}}
\newcommand{\wt}{\mathrm{wt}}
\newcommand{\slwt}{\wt_{\mathrm{sl}}}
\newcommand{\Wb}{\overline{W}}
\newcommand{\La}{\Lambda}
\newcommand{\Lab}{\overline{\Lambda}}
\newcommand{\la}{\lambda}
\newcommand{\lab}{\overline{\lambda}}
\newcommand{\lh}{\widehat{\la}}
\newcommand{\alb}{\overline{\alpha}}
\newcommand{\Qh}{\widehat Q}
\newcommand{\Rh}{\widehat{R}}
\newcommand{\bh}{\widehat{b}}
\newcommand{\BB}{\mathcal{B}}
\newcommand{\VV}{\mathcal{V}}
\newcommand{\dotnu}{{\nu^{\cdot}}}
\newcommand{\C}{{\mathbb C}}
\newcommand{\Z}{{\mathbb Z}}
\newcommand{\Q}{{\mathbb Q}}
\newcommand{\N}{\mathbb Z_+}
\newcommand{\inner}[2]{\langle #1,#2\rangle}
\newcommand{\K}{K}
\newcommand{\LL}{\mathcal{L}}
\newcommand{\PPP}{\mathbb{P}}
\newcommand{\QQQ}{\mathbb{Q}}
\newcommand{\et}{\widetilde{e}}
\newcommand{\ft}{\widetilde{f}}
\newcommand{\st}{\widetilde{r}}
\newcommand{\bb}{\overline{b}}
\newcommand{\Qb}{\overline{Q}}
\newcommand{\tK}{\widetilde{K}}
\newcommand{\SLNHAT}{\widehat{SL}_n}
\newcommand{\slnhat}{\widehat{sl}_n}
\newcommand{\rows}{\mathrm{rows}}
\newcommand{\PSE}{\mathrm{P}_\searrow}
\newcommand{\shape}{\mathrm{shape}}
\newcommand{\LRT}{\mathrm{LRT}}
\newcommand{\CST}{\mathrm{CST}}
\newcommand{\pr}{\mathrm{pr}}
	
\section{Introduction}

In \cite[Theorem 5.2]{KMOTU} it was shown that the characters of
certain level one Demazure modules of type $A^{(1)}_{n-1}$, when decomposed
as linear combinations of irreducible characters of type $A_{n-1}$,
have coefficients given by Kostka-Foulkes polynomials in the variable
$q=e^{-\delta}$ where $\delta$ is the null root.  The key steps in the
proof are that
\begin{enumerate}
\item The Demazure crystals are isomorphic to tensor products
of classical $\slnhat$ crystals indexed by fundamental weights \cite{KMOU}.
\item The generating function over these crystals by weight
and energy function is equal to the generating function over
column-strict Young tableaux by weight and charge \cite{NY}.
\item The Kostka-Foulkes polynomial is the coefficient of an irreducible
$sl_n$-character in the tableau generating function \cite{Bu} \cite{LS1}.
\end{enumerate}

The main result of this paper is that for a Demazure module
of arbitrary level whose lowest weight is a multiple of one of those
in \cite{KMOTU}, the corresponding coefficient polynomial is the
Poincar\'e polynomial of an isotypic component of the coordinate ring
of the closure of the conjugacy class of a nilpotent matrix.
The Poincar\'e polynomial is the $q$-analogue of the multiplicity
of an irreducible $gl(n)$-module in a tensor product
of irreducible $gl(n)$-modules in which each factor has highest weight
given by a rectangular (that is, a multiple of a fundamental) weight.
These polynomials, which possess many properties generalizing those
of the Kostka-Foulkes, have been studied extensively using
algebro-geometric and combinatorial methods \cite{KKSW} \cite{KS}
\cite{S1} \cite{S2} \cite{SW1} \cite{SW2}.

The connection between the Demazure modules and the nilpotent adjoint
orbit closures can be explained as follows.  Let $X_\mu$ be the
Zariski closure of the conjugacy class of the nilpotent Jordan matrix
with block sizes given by the transpose partition $\mu^t$ of $\mu$,
that is,
\begin{equation*}
  X_\mu = \{ A \in \mathrm{gl}_n(\C) \mid 
	\mathrm{dim\,\, ker\,\, }A^i \ge \mu_1+\dots+\mu_i
\text{ for all $i$ } \}
\end{equation*}
Lusztig gave an embedding of the variety $X_\mu$
as an open dense subset of a $P$-stable Schubert variety
$Y_\mu$ in $\SLNHAT/P$, where $P\cong SL_n$ is the 
parabolic subgroup given by ``omitting the reflection $r_0$" \cite{Lu}.
The desired level $l$ Demazure module, viewed as an $sl_n$-module,
is isomorphic to the dual of the space 
of global sections $H^0(Y_\mu,\LL_{\La_0}^{\otimes l})$,
where $\LL_{\La_0}$ is the restriction to $Y_\mu$ of the homogeneous
line bundle on $\SLNHAT/P$ affording the fundamental weight $\La_0$.

The proof of the main result entails generalizations of the
three steps in the proof of \cite[Theorem 5.2]{KMOTU}.  First,
the methods of \cite{KMOTU} may used to
show that the Demazure crystal is isomorphic to a classical $\slnhat$
crystal $B$ that is a tensor product of crystals of the form $B^{k,l}$
(notation as in \cite{KMN2}).
We call $B^{k,l}$ a rectangular crystal since it is indexed by
the weight $l\La_k$ that corresponds to the rectangular partition
with $k$ rows and $l$ columns.

Second, it is shown that the crystal $B$ is indexed by
sequences of Young tableaux of rectangular shape equipped
with a generalized charge map.  In particular, 
we give explicit descriptions in terms of tableaux and the
Robinson-Schensted-Knuth correspondence, of
\begin{enumerate}
\item[$\bullet$] The zero-th crystal raising operator $\et_0$ acting on $B$,
which involves the generalized cyclage operators of \cite{S1}
on LR tableaux and a promotion operator on column-strict tableaux.
\item[$\bullet$] The combinatorial $R$-matrices
on a tensor product of the form $B^{k_1,l_1}\otimes B^{k_2,l_2}$,
which are given by a combination of the generalized automorphism
of conjugation \cite{S1} and the energy function.
\item[$\bullet$] The energy function, which equals the generalized
charge of \cite{KS} \cite{S1}.
\end{enumerate}
Moreover, it is shown that every generalized cocyclage relation
\cite{S1} on LR tableaux may be realized by $\et_0$.  The
formula for the corresponding change in the energy function by
$\et_0$ was known \cite{NY} in the case that all rectangles
are single rows or all are single columns.

Third, it must be shown that the tableau formula coincides with the
Poincar\'e polynomial.  This was accomplished in
\cite{S1}, where it is shown that the tableau formula satisfies
a defining recurrence of Weyman \cite{SW2} \cite{W} for the Poincar\'e
polynomials that is closely related to Morris' recurrence for Kostka-Foulkes
polynomials \cite{Mo}.  

As an application of the formula for the energy function, we
give a very simple proof of a monotonicity property
for the Poincar\'e polynomials (conjectured by A. N. Kirillov)
that extends the monotonicity property of the Kostka-Foulkes polynomials
that was proved by Han \cite{Ha}.

Thanks to M. Okado for pointing out the reference \cite{ScWa}
which has considerable overlap with this paper and \cite{S1}.

\section{Notation and statement of main result}

\subsection{Quantized universal enveloping algebras}
For this paper we only require the following three algebras:
$U_q(sl_n)\subset U'_q(\slnhat)\subset U_q(\slnhat)$.

Let us recall some definitions for quantized universal enveloping algebras
taken from \cite{Kac} and \cite{KMN}.  Consider the following data:
a finitely generated $\Z$-module $P$ (weight lattice),
a set $I$ (index set for Dynkin diagram), elements
$\{\alpha_i\mid i\in I\}$ (basic roots) and 
$\{\ h_i\in P^*=\Hom_{\Z}(P,\Z) \mid i\in I\}$ (basic coroots)
such that $(\inner{h_i}{\alpha_j})_{i,j\in I}$ is a generalized Cartan
matrix, and a symmetric form $(\cdot,\cdot): P \times P \rightarrow \Q$
such that $(\alpha_i,\alpha_i)\in \Z$ is positive,
$\inner{h_i}{\la}=2(\alpha_i,\la)/(\alpha_i,\alpha_i)$ for
$i\in I$ and $\la \in \Q\otimes P$.

This given, let $U_q(\slnhat)$ be the quantized universal
enveloping algebra, the $\Q(q)$-algebra with generators
$\{e_i,f_i\mid i\in I\}$ and
$\{q^h\mid h\in P^*\}$ and relations as in \cite[Section 2]{KMN}.

For $U_q(\slnhat)$, let $I=\{0,1,2\dots,n-1\}$
be the index set for the Dynkin diagram, $(a_{ij})_{i,j\in I}$ the
Cartan matrix of type $A^{(1)}_{n-1}$, $P$ the free $\Z$-module with basis
$\{ \La_i \mid i\in I\} \cup \{\delta\}$ (fundamental weights)
and let $P^*$ have dual basis $\{h_i\mid i\in I\}\cup\{d\}$. 
Define the elements $\{\alpha_i\in P \mid i\in I\}$ by
\begin{equation*}
  \alpha_i = \delta_{0i}\delta + \sum_{j\in I} a_{ij} \La_j,
\end{equation*}
so that $(\inner{h_i}{\alpha_j})_{i,j\in I}$ is the Cartan
matrix of type $A^{(1)}_{n-1}$ and $\delta=\sum_{i=0}^{n-1} \alpha_i$.
Define the symmetric $\Q$-valued form
$(\cdot,\cdot)$ by $(\alpha_i,\alpha_j)=a_{ij}$ for
$i,j\in I$, $(\alpha_i,\La_0)=\delta_{i0}$ for $i\in I$,
and $(\La_0,\La_0)=0$.  Then the quantized universal enveloping
algebra for this data is $U_q(\slnhat)$.
Let $P^+ = \N\delta \oplus \bigoplus_{i\in I} \N \La_i$
be the dominant weights.  For $\La\in P^+$ let $\VV(\La)$
be the irreducible integrable highest weight $U_q(\slnhat)$-module
of highest weight $\La$, $\BB(\La)$ its crystal graph,
and $u_\La\in\BB(\La)$ the highest weight vector.

For $U'_q(\slnhat)$, let $I$ and $(a_{ij})$ be as above,
but instead of $P$ use the ``classical weight lattice"
$\PCL=P/\Z\delta=\bigoplus_{i\in I} \Z \La_i$ (where
by abuse of notation the image of $\La_i$ in $\PCL$ is also denoted
$\La_i$).  The basic coroots $\{h_i\mid i\in I\}$ form a $\Z$-basis
of $\PCL^*$.  The basic roots are $\{\alb_i\mid i\in I\}$
where $\alb_i$ denotes the image of $\alpha_i$ in $\PCL$ for $i\in I$.
Note that the basic roots are linearly dependent.
The pairing and symmetric form are induced by those above.
The algebra for this data is denoted $U'_q(\slnhat)$ and may be viewed
as a subalgebra of $U_q(\slnhat)$ since its generators are a subset
of those of $U_q(\slnhat)$ and its relations map to relations.
Let $\PCL^+ = \bigoplus_{i\in I} \N \La_i$.
For $\La\in P^+$, the $U_q(\slnhat)$-module
$\VV(\La)$ is a $U'_q(\slnhat)$-module by restriction
and weights are taken modulo $\delta$.

For $U_q(sl_n)$, let $J=\{1,2,\dots,n-1\}\subset I$ be the index
set for the Dynkin diagram, $(a_{ij})_{i,j\in J}$ the restriction
of the above Cartan matrix to $J\times J$, and $\PCLB=\PCL/\Z \La_0$
the weight lattice.  The basic coroots $\{h_i\in\PCLB^*\mid i\in J\}$.
form a $\Z$-basis of $\PCLB^*=\bigoplus_{i\in J}\Z h_i$.  The basic
roots are $\{\alb_i\in\PCLB\mid i\in J\}$.
The algebra for this data is $U_q(sl_n)$, which can be viewed as a subalgebra
of $U'_q(\slnhat)$.  Let $\PCLB^+=\bigoplus_{i\in J} \N \La_i$ be the
dominant integral weights.  For $\la\in\PCLB^+$ let $V^\la$ be
the irreducible $U_q(sl_n)$-module of highest weight $\la$,
and $B^\la$ its crystal graph.

Denote by $W$ the Weyl group of the algebra $\slnhat$ and
by $\Wb$ that of $sl_n$.  $W$ is the subgroup of automorphisms of $P$
generated by the simple reflections $\{ r_i \mid i\in I\}$ where
\begin{equation*}
  r_i(\la) = \la - \inner{h_i}{\la} \,\alpha_i.
\end{equation*}
Let $\Qb=\bigoplus_{i\in J} \Z \alb_i \subset\PCL$.
$W$ acts faithfully on the affine subspace
$X = \La_0 + \Qb \subset \PCL$.  For $\la\in\Qb$ write $\tau_\la$ for
the map $X\rightarrow X$ given by translation by $\la$.
Let $\theta =\sum_{i\in J} \alb_i \in \PCL$.
Then the action of $r_i$ on $X$ for $i \in J$ is given by
\begin{equation*}
  r_i(\lab) = \lab - \inner{h_i}{\lab}\, \alb_i,
\end{equation*}
for $\lab\in X$ and $r_0 = \tau_{\theta} r_\theta$, where
\begin{equation*}
  r_\theta = r_1 \dots r_{n-2} r_{n-1} r_{n-2} \dots r_1
\end{equation*}
is the reflection through the hyperplane orthogonal to $\theta$.
Then $W \cong \Wb \ltimes \Qb$ where the element
$\la\in \Qb$ acts by $\tau_\la$.

For $\La\in P^+$ and $w\in W$ the Demazure module of lowest weight
$w\La$ is defined by $\VV_w(\La) = U_q(\mathfrak{b}) v_{w \La}$
where $v_{w \La}$ is a generator of the (one dimensional) extremal weight
space in $\VV(\La)$ of weight $w\La$ and $U_q(\mathfrak{b})$ is the
subalgebra of $U_q(\slnhat)$ generated by the $e_i$ and $h\in P^*$.

\subsection{Main result}

Let $\mu$ be a partition of $n$.  The coordinate ring $\C[X_\mu]$ of the
nilpotent adjoint orbit closure $X_\mu$ has a graded $sl_n$-action
induced by matrix conjugation on $X_\mu$.  For $\la\in\PCLB^+$,
define the Poincar\'e polynomial of the $\la$-th isotypic
component of $\C[X_\mu]$, by
\begin{equation*}
  \K_{\la;\mu}(q) = \sum_{d\ge0} q^d \dim \Hom_{sl_n}(V^\la,\C[X_\mu]_d)
\end{equation*}
where $\C[X_\mu]_d$ is the homogeneous component of degree $d$.

Partitions with at most $n$ parts are projected to dominant
integral weights of $sl_n$ by $\mu\mapsto \slwt(\mu)$, where
\begin{equation*}
 \slwt(\mu)=\sum_{i=1}^{n-1} (\mu_i-\mu_{i+1})\La_i\in\PCLB^+.
\end{equation*}

\begin{rem} Warning: this is not the Kostka polynomial, but a generalization;
see \cite{SW2}.  In the special case that
$\la$ is a partition of $n$ with at most $n$ parts
then $\K_{\slwt{\la};\mu}(q) = \tK_{\la^t\mu}(q)$ which is a
renormalization of the Kostka-Foulkes polynomial with indices
$\la^t$ and $\mu$.
\end{rem}

\begin{thm} \label{main result} Let $l$ be a positive integer,
$\mu$ a partition of $n$, and $w_\mu\in W$ the translation by
the antidominant weight $-w_0 \slwt(\mu^t)\in \PCLB$,
where $w_0$ is the longest element of $\Wb$.  Then
\begin{equation*}
  e^{-l \La_0}\, \ch \VV_{w_\mu}(l \La_0) = 
   \sum_\la \K_{\slwt(\la);\mu}(q)\, \ch V^{\slwt(\la)}
\end{equation*}
where $\la$ runs over the partitions of the multiple
$\,l n$ of $n$ with at most $n$ parts.
\end{thm}

\section{Crystal structure on tensor products of rectangles}

The goal of this section is to give explicit descriptions of
the classical $\slnhat$ crystal structure on
tensor products of rectangular crystals and their energy functions.
This is accomplished by translating the theory of $sl_n$ crystals
and classical $\slnhat$ crystals in \cite{KMN} \cite{KMN2} \cite{KN}
\cite{NY} into the language of Young tableaux and the
Robinson-Schensted-Knuth (RSK) correspondence.

\subsection{Crystals}
This section reviews the definition of a weighted crystal \cite{KMN}
and gives the convention used here for the tensor product of crystals.

A $P$-weighted $I$-crystal is a a weighted $I$-colored directed graph $B$,
that is, a set equipped with a weight function $\wt:B\rightarrow P$
and directed edges colored by the set $I$, satisfying the following
properties.
\begin{enumerate}
\item[(C1)] There are no multiple edges; that is, for each $i\in I$
and $b,b'\in B$ there is at most one edge colored $i$ from $b$ to $b'$.
\end{enumerate}
If such an edge exists, this is denoted $b'=\ft_i(b)$ or equivalently
$b=\et_i(b')$, by abuse of the notation of a function $B\rightarrow B$.
It is said that $\ft_i(b)$ is defined or equivalently that $\et_i(b')$ is
defined, if the edge exists.
\begin{equation*}
\begin{split}
  \phi_i(b) &= \max \{ m\in\N \mid \text{$\ft_i^m(b)$ is defined} \} \\
  \epsilon_i(b) &= \max \{ m\in\N \mid \text{$\et_i^m(b)$ is defined.} \}
\end{split}
\end{equation*}
\begin{enumerate}
\item[(C2)] If $\ft_i(b)$ is defined then $\wt(\ft_i(b))=\wt(b)-\alpha_i$.
Equivalently, $\wt(\et_i(b))=\wt(b)+\alpha_i$.
\item[(C3)] $\inner{h_i}{\wt(b)}=\phi_i(b)-\epsilon_i(b)$.
\end{enumerate}

If $B_j$ is a $P$-weighted $I$-crystal for $1\le j\le m$,
the Cartesian product $B_m \times \dots \times B_1$ can be given
a crystal structure as follows; this crystal is denoted
$B=B_m\otimes \dots\otimes B_1$.  The convention used here is opposite
that in much of the literature  but is convenient for the tableau
combinatorics used later.  Let $b_j\in B_j$ and
$b=b_m\otimes \dots\otimes b_1\in B$.  The weight function on $B$ is
given by
\begin{equation*}
  \wt(b) = \sum_{j=1}^m \wt(b_j).
\end{equation*}
The root operators $\ft_i$ and the functions $\phi_i$ are defined by the
``signature rule".
Given $b\in B$ and $i\in I$, construct a biword (sequence of pairs of
letters) consisting of $\phi_i(b_j)$ copies of the
biletter $\binom{j}{-}$ and $\epsilon_i(b_j)$ copies of the
biletter $\binom{j}{+}$ for all $j$, sorted in weakly increasing order
by the order $\binom{j}{\pm}<\binom{j'}{\pm}$ if $j>j'$ and
$\binom{j}{-}<\binom{j}{+}$.  This biword is now repeatedly reduced by
removing adjacent biletters whose lower letters are $+-$ in that order.
If $+$ and $-$ are viewed as left and right parentheses then this
removes matching pairs of parentheses.  At the end one obtains a biword
whose lower word has the form $-^s+^t$.  Then define
$\phi_i(b)=s$ and $\epsilon_i(b)=t$.  If $s>0$
(resp. $t>0$) let $j_-$ (resp. $j_+$) be the upper
letter corresponding to the rightmost $-$ (resp. leftmost $+$)
in the reduced biword, and define
\begin{equation*}
  \ft_i(b) = b_m \otimes \dots \otimes b_{1+j_-} \otimes
  	\ft_i(b_{j_-}) \otimes b_{-1+j_-} \otimes \dots \otimes b_1,
\end{equation*}
respectively,
\begin{equation*}
  \et_i(b) = b_m \otimes \dots \otimes b_{1+j_+} \otimes
  	\et_i(b_{j_+}) \otimes b_{-1+j_+} \otimes \dots \otimes b_1.
\end{equation*}
A morphism $g:B\rightarrow B'$ of $P$-weighted $I$-crystals is a map $g$
that preserves weights and satisfies $g(\ft_i(b))=\ft_i(g(b))$
for all $i\in I$ and $b\in B$, that is, if $\ft_i(b)$ is defined
then $\ft_i(g(b))$ is, and the above equality holds.

It is easily verified that the $P$-weighted $I$-crystals
form a tensor category.

We only require the following kinds of crystals.
\begin{enumerate}
\item The crystal graphs of integrable $U_q(\slnhat)$-modules are
$P$-weighted $I$-crystals and are called $\slnhat$-crystals.
\item The crystal graphs of $U'_q(\slnhat)$-modules that are either
integrable or are finite-dimensional and have a weight space
decomposition, are $\PCL$-weighted $I$-crystals
and are called classical $\slnhat$-crystals.
\item The crystal graphs of integrable $U_q(sl_n)$-modules are
$\PCLB$-weighted $J$-crystals and are called $sl_n$-crystals.
\end{enumerate}

\subsection{Crystal reflection operator and Weyl group action}
Let $B$ be an $sl_n$ crystal and $i\in J$.
Write $p=\phi_i(b)-\epsilon_i(b)$.  Define
\begin{equation*}
  \st_i(b) = \begin{cases}
    \ft_i^p(b) & \text{if $p>0$} \\
    b & \text{if $p=0$} \\
    \et_i^{\,-p}(b) & \text{if $p<0$}
  \end{cases}
\end{equation*}
The Weyl group $\Wb$ acts on $B$ by $r_i(b)=\st_i(b)$
for $i\in J$.  It is obvious that $\st_i$ is an involution and that
$\st_i$ and $\st_j$ commute for $|i-j|>1$, but
not at all obvious that the $\st_i$ satisfy the other braid relation.
A combinatorial proof of this fact is indicated in
\cite{LS2} for the action of $\Wb$ on an irreducible $sl_n$ crystal.

\subsection{Irreducible $sl_n$ crystals}
\label{sl_n sec}

Let $\la=(\la_1\ge\la_2\ge\dots\ge\la_n)\in\N^n$
be a partition of length at most $n$.
Let $V^\la$ be the irreducible $U_q(sl_n)$-module of highest
weight $\slwt(\la)$ and $B^\la$ its crystal.
In \cite{KN} the structure of the $sl_n$ crystal $B^\la$
is determined explicitly.  The crystal $B^\la$ 
may be indexed by the set $\CST(\la)$ of column-strict tableaux of
shape $\la$ with entries in the set $[n]=\{1,2,\dots,n\}$.
The combinatorial construction yielding the action of the
crystal operators $\et_i$ and $\ft_i$ for $1\le i\le n-1$ on
tableaux was already known.  In a 1938 paper, in the course of
proving the Littlewood-Richardson rule, G. de B. Robinson
gave a form of the Robinson-Schensted-Knuth (RSK) correspondence
which is defined by giving the value of the map on $sl_n$-highest
weight vectors and then extending it by via canonical sequences of raising
operators $\ft_i$ \cite[Section 5]{R}; see also \cite[I.9]{Mac} where
Robinson's proof is cleaned up.

Suppose first that $\la=(1)$ so $B^\la$ is the crystal of the
defining representation of $sl_n$.  This crystal is indexed
by the set $[n]=\{1,2,\dots,n\}$ and $\ft_i(j)$ is defined if and
only if $j=i$ and in that case $\ft_i(i)=i+1$.

Next consider the tensor product $(B^{(1)})^{\otimes m}$.
It may be indexed by words $u=u_m\dots u_1$ of length $m$
in the alphabet $[n]$, where $u_j\in [n]$.
Its $sl_n$ crystal structure is defined by the signature rule.
This case of the signature rule is given in \cite{LS2}.

Now it is necessary to introduce notation for Young tableaux.

Some definitions are required.
The Ferrers diagram $D(\la)$ is the set of pairs of integers
$D(\la) = \{(i,j)\in\N^2 \mid 1\le j\le \la_i \}$.
A skew shape $D=\la/\mu$ is the set difference of the Ferrers diagrams
$D(\la)$ and $D(\mu)$ of the partitions $\la$ and $\mu$.
If $D$ and $E$ are skew shapes such that $D$ has $c$ columns
and $E$ has $r$ rows, then define the skew shape
\begin{equation*}
  D \otimes E = \{(i+r,j)\mid (i,j)\in D\} \cup
  	\{(i,j+c)\mid (i,j)\in E\}.
\end{equation*}
In other words, $D\otimes E$ is the union of a translate of $D$
located to the southwest of a translate of $E$.

A tableau of (skew) shape $D$ is a function
$T:D\rightarrow[n]=\{1,2,\dots,n\}$, and is depicted as a partial
matrix whose $(i,j)$-th entry is $T(i,j)$ for all $(i,j)\in D$.
Denote by $\shape(T)$ the domain of $T$.
The tableau $T$ is said to be column-strict if $T(i,j)\le T(i,j+1)$
for all $(i,j),(i,j+1)\in D$ and $T(i,j)<T(i+1,j)$ for all
$(i,j),(i+1,j)\in D$.
Let $\CST(D)$ be the set of column-strict tableaux of shape $D$.
The content of a tableau is the sequence
\begin{equation*}
\content(T)=(m_1(T),\dots,m_n(T))
\end{equation*}
where $m_i(T)$ is the number of occurrences of the letter $i$ in $T$.
Let $\CST(D,\gamma)$ denote the set of column-strict tableaux of shape
$D$ and content $\gamma$.
The (row-reading) word of the tableau $T$ is the word given by
$\word(T)=\cdots v^2 v^1$ where $v^i$ is the word obtained
by reading the $i$-th row of $T$ from left to right.
Say that the word $u$ fits the skew shape $D$ if there is a
column-strict tableau (necessarily unique) whose row-reading word is $u$.

\begin{rem} \label{skew word} Let $D$ be a skew shape, $u$ a word
in the alphabet $[n]$ and $1\le i\le n-1$.  It is well-known and easy
to verify that if $\et_i(u)$ is defined, then $u$ fits $D$
if and only if $\et_i(u)$ does.
This given, if $T$ is a column-strict tableau of shape $D$ and
$\et_i(\word(T))$ is defined, then let $\et_i(T)$ be the unique
column-strict tableau of shape $D$ such that
$\word(\et_i(T))=\et_i(\word(T))$.  The same can be done for $\ft_i$.
\end{rem}

Thus the set $\CST(D)$ is an $sl_n$ crystal; call it $B^D$.
When $D=D(\la)$ this is the crystal $B^\la$.

\subsection{Tensor products of irreducible $sl_n$ crystals and RSK}

Let $D$ be a skew diagram and $B^D$ the $sl_n$ crystal defined in the
previous section.  The RSK correspondence yields a combinatorial
decomposition of $B$ into irreducible $sl_n$ crystals.
The RSK map can be applied to
tensor products of irreducible $sl_n$ crystals.
The goal of this section is to review a well-known parametrizing
set for the multiplicity space of such a tensor product, by what
we shall call Littlewood-Richardson (LR) tableaux.

Let $\eta=(\eta_1,\eta_2,\dots,\eta_m)$ be a sequence of
positive integers summing to $n$ and $R=(R_1,R_2,\dots,R_m)$ a
sequence of partitions such that $R_i$ has $\eta_i$ parts, some of
which may be zero.  Let $A_1$ be the first $\eta_1$ numbers in $[n]$,
$A_2$ the next $\eta_2$, and so on.  Recall the skew shape
$R_m\otimes\dots \otimes R_1$, embedded in the plane so that
$A_i$ gives the set of row indices for $R_i$.  Let $\gamma=(\gamma_1,\dots,
\gamma_n)\in\N^n)$ where $\gamma_i$ is the length of the $i$-th row of
$R_m\otimes\dots\otimes R_1$, that is, $\gamma$ is obtained by juxtaposing
the partitions $R_1$ through $R_m$.  The tensor product crystal
may be viewed as a skew crystal:
\begin{equation*}
  B^R := B^{R_m\otimes\dots\otimes R_1} \cong B^{R_m} \otimes \dots
   	\otimes B^{R_1}.
\end{equation*}
Let $b\in B^R$ and write $b = b_m\otimes\dots\otimes b_1$
where $b_i\in\CST(R_i)$.  Let $v_r$ for the weakly increasing
word (of length $\gamma_r$) comprising the $r$-th row of $b$, for
$1\le r\le n$.  The word of $b$ regarded as a skew column-strict tableau,
is given by $\word(b)=v_n \cdots v_1$.

Recall Knuth's equivalence on words \cite{Knuth}.  Say that a skew
shape is normal (resp. antinormal) if it has a unique northwest
(resp. southeast) corner cell \cite{H}.  A normal skew shape is merely a
translation in the plane of a partition shape, and an antinormal shape
is the 180-degree rotation of a normal shape.  For any word $v$, there is
a unique (up to translation) column-strict tableau $P(v)$ of normal shape 
such that $v$ is Knuth equivalent to the word of $P(v)$.
There is also a unique (up to translation) skew column-strict tableau
$\PSE(v)$ of antinormal shape such that $v$ is Knuth
equivalent to the word of $\PSE(v)$.

The tableau $P(v)$ may be computed by Schensted's column-insertion
algorithm \cite{Schen}.  For a subinterval $A\subset[n]$ and a (skew)
column-strict tableau $T$, define $T|_A$ to be the skew column-strict tableau
obtained by restricting $T$ to $A$, that is, removing from $T$ the letters
that are not in $A$.  Define the pair of column-strict tableaux
$(\PPP(b),\QQQ(b))$ by
\begin{equation*}
\begin{split}
  \PPP(b) &= P(\word(b)) \\
  \shape(\QQQ(b)|_{[r]}) &= \shape(P(v_r \dots v_1))
	\qquad\text{for all $0\le r\le n$}
\end{split}
\end{equation*}
By definition $\shape(\PPP(b))=\shape(\QQQ(b))$.
It is easy to show that $\QQQ(b)$ is column-strict and of content $\gamma$.
This gives an embedding
\begin{equation} \label{RSK}
\begin{split}
  B^R &\hookrightarrow \bigcup_\la \CST(\la)\times \CST(\la,\gamma) \\
  b &\mapsto (\PPP(b),\QQQ(b))
\end{split}
\end{equation}
It is well-known that this is a map of $sl_n$ crystals.
That is, if $g$ is any of $\et_i$, $\ft_i$, or $\st_i$ for $i\in J$, then
\begin{equation} \label{crystal RSK}
\begin{split}
  \PPP(g(b))&=g(\PPP(b)) \\
  \QQQ(g(b))&=\QQQ(b)
\end{split}
\end{equation}
For the case $g=\st_i$ this fact is in \cite{LS2}.

Let us describe the image of the map \eqref{RSK}.
For the partition $\la=(\la_1,\dots,\la_n)$ and
a permutation $w$ in the symmetric group $S_n$,
the key tableau $\key(w \la)$ of content $w \la$, is the
unique column-strict tableau of shape $\la$ and content $w \la$.
In the above notation for the sequence of partitions $R$,
for $1\le j\le m$ let $Y_j=\key(R_j)$ in the alphabet $A_j$.

Say that a word $u$ in the alphabet $[n]$ is $R$-LR (short for
$R$-Littlewood-Richardson) if $P(u|_{A_j})=Y_j$ for all $1\le j\le m$,
where $u|_{A_j}$ is the restriction of the word $u$ to the
subalphabet $A_j\subset[n]$.  Say that a (possibly skew) column-strict
tableau is $R$-LR if its row-reading word is.  Denote by
$\LRT(\la;R)$ the $R$-LR tableaux of partition shape $\la$
and $\LRT(R)=\bigcup_\la \LRT(\la;R)$.

The following theorem is essentially a special case of
\cite[Theorem 1]{W}  which is a strong version
of the classical rule of Littlewood and Richardson \cite{LR}.

\begin{thm} The map \eqref{RSK} gives a bijection
\begin{equation} \label{R RSK}
  B^R \cong \bigcup_\la \CST(\la)\times \LRT(\la;R).
\end{equation}
\end{thm}

\subsection{Kostka crystals}
\label{Kostka sec}

In the case that $m=n$ and $R_i=(\gamma_i)$ for
$\gamma\in\N^n$, we call $B^R$ a Kostka crystal.  The set
$\LRT(\la;R)$ is merely the set $\CST(\la,\gamma)$.

Let $R=(R_1,\dots,R_m)$ be a sequence of rectangles as usual.
Define the Kostka crystal $B^{\rows(R)}$ by the
sequence of one row partitions $\rows(R)$ whose $r$-th
partition is given by $(\gamma_r)$ where $\gamma_r$ is
the length of the $r$-th row of the skew shape
$R_m\otimes \dots \otimes R_1$.  Letting $b\in B=B^R$ and
$v_r$ the $r$-th row of $b$, there is the obvious
$sl_n$ crystal embedding
\begin{equation} \label{embed into Kostka}
\begin{split}
  B^R &\hookrightarrow B^{\rows(R)} \\
  b &\mapsto v_n \otimes \dots\otimes v_1.
\end{split}
\end{equation}

In fact, the RSK correspondence \eqref{R RSK} may be defined using
the commutativity of the diagram in \eqref{crystal RSK} for
$g=\ft_i$ and all $i\in J$, and giving its values on
the $sl_n$-highest weight elements in $B$ \cite{R}.
Suppose $b\in B$ is such that $\word(b)$ is of $sl_n$-highest weight,
that is, $\epsilon_i(b)=0$ for all $i\in J$.
Such words are said to possess the lattice property.  In this case
$\content(\word(b))$ must be a partition, say $\la$, and
$\PPP(b)=P(\word(b))=\key(\la)$.
For the recording tableau, write $\word(b)=v_n\otimes \dots\otimes v_1$
where $v_r$ is the $r$-th row of $b$ viewed as a tableau of the skew
shape $R_m\otimes\dots\otimes R_1$.  Then
$\QQQ(b)$ is the column-strict tableau of shape equal to
$\content(b)$, whose $i$-th row contains $m$ copies of the letter $j$
if and only if the word $v_j$ contains $m$ copies of the letter $i$,
for all $i$ and $j$.
In particular, if elements $b,b'\in B$ admit
the same sequences of raising operators, then $\PPP(b)=\PPP(b')$.

\subsection{Rectangle-switching bijections}
From now on we consider only crystals $B^R$ where
$R_j=(\mu_j^{\eta_j})$ is the partition with $\eta_j$ rows and
$\mu_j$ columns for $1\le j\le m$.

Consider the case $m=2$, with $R=(R_1,R_2)$, $r_1 R = (R_2,R_1)$,
$B=B^R=B^{R_2}\otimes B^{R_1}$ and $r_1 B=B^{r_1 R}=B^{R_1}\otimes B^{R_2}$.
Since the $U_q(sl_n)$-module
$V^{R_2}\otimes V^{R_1}\cong V^{R_1}\otimes V^{R_2}$ is multiplicity-free,
it follows that there is a unique $sl_n$ crystal isomorphism
\begin{equation*}
 \sigma_{R_2,R_1}: B^{R_2}\otimes B^{R_1}\cong B^{R_1}\otimes B^{R_2}.
\end{equation*}
It is defined explicitly in terms of the RSK correspondence as follows.
By the above multiplicity-freeness, for any partition $\la$,
\begin{equation*}
|\LRT(\la;(R_1,R_2))| = \LRT(\la;(R_2,R_1))| \le 1.
\end{equation*}
Thus there is a unique bijection
\begin{equation*}
\tau=\tau_{\la;(R_2,R_1)}:\LRT(\la;(R_1,R_2))\rightarrow \LRT(\la;(R_2,R_1)).
\end{equation*}
For later use, extend $\tau$ to a bijection from the set of $(R_1,R_2)$-LR
words to the set of $(R_2,R_1)$-LR words by
\begin{equation} \label{tau-prime}
\begin{split}
  P(\tau(u))&=\tau(P(u)) \\
  Q(\tau(u))&=Q(u)
\end{split}
\end{equation}
where $Q(u)$ is the standard column-insertion recording tableau
(that is, $Q(u)=\QQQ(u)$ with $u$ regarded as a word in the
tensor product $(B^{(1)})^{\otimes N}$ where $N$ is the length of $u$).

$\tau$ is the rectangular generalization of an automorphism of conjugation.
$\sigma$ is defined by the commutative diagram
\begin{equation*}
\begin{CD}
  B^{R_1\otimes R_2} @>\mathrm{RSK}>> \bigcup_\la \CST(\la)\times 
	\LRT(\la;(R_1,R_2)) \\
   @V{\sigma}VV @VV{\bigcup_\la \id \times \tau}V  \\
  B^{R_2\otimes R_1} @>\mathrm{RSK}>> \bigcup_\la \CST(\la)\times 
	\LRT(\la;(R_2,R_1))
\end{CD}
\end{equation*}
In other words, for all $b\in B^{R_2} \otimes B^{R_1}$,
\begin{equation} \label{rect switch}
\begin{split}
\PPP(\sigma(b))&=\PPP(b) \\
\QQQ(\sigma(b))&=\tau(\QQQ(b)).
\end{split}
\end{equation}

Now let $R=(R_1,\dots,R_m)$ and $B=B^R$ the tensor product of
rectangular crystals.  Let $w\in S_m$ be a permutation in the symmetric
group on $m$ letters.  Write $w B=B^{w R}$ where $w R$ is the
sequence of rectangles
$w R=R_{w^{-1}(m)} \otimes \dots \otimes R_{w^{-1}(m)}$.

Write $\sigma_{R_i,R_j}$ for the action of the above $sl_n$ crystal
isomorphism at consecutive tensor positions in
$w B = \dots \otimes B^{R_i} \otimes B^{R_j} \otimes \dots$.
Then the isomorphisms $\sigma_{R_i,R_j}$ satisfy a Yang-Baxter identity
\begin{equation*}
\begin{split}
      &(\sigma_{B_j,B_k} \otimes \id)
\circ (\id \otimes \sigma_{B_i,B_k})
\circ (\sigma_{B_i,B_j} \otimes \id) \\
=\,\,
      &(\id \otimes \sigma_{B_i,B_j})
\circ (\sigma_{B_i,B_k} \otimes \id)
\circ (\id \otimes \sigma_{B_j,B_k})
\end{split}
\end{equation*}
This is a consequence of the corresponding difficult identity for
bijections $\tau_{R_i,R_j}$ on recording tableaux,
defined and conjectured in \cite{KS} and proven in \cite[Theorem 9 (A5)]{S1}.
By composing maps of the form $\sigma_{R_i,R_j}$, it is possible
to well-define $sl_n$ crystal isomorphisms
\begin{equation} \label{sigma perm}
\sigma_{R,wR}: B \rightarrow w B
\end{equation}
such that $\sigma_{R,vw R}=\sigma_{w R,v w R}\circ \sigma_{R,w R}$
for $v,w\in S_m$.  These bijections satisfy
\begin{equation} \label{sigma tau perm}
\begin{split}
  \PPP(\sigma_{R,wR}(b)) &= \PPP(b) \\
  \QQQ(\sigma_{R,wR}(b)) &= \tau_{R,wR}(\QQQ(b))
\end{split}
\end{equation}
where $\tau_{R,wR}:\LRT(R)\cong\LRT(wR)$ is the
shape-preserving bijection defined in \cite{S1}.

\begin{rem} \label{coreflection}
Suppose $B=B^R$ is a Kostka crystal,
$b\in B$, and $b=b_n \otimes \dots \otimes b_1$ with
with $R_j=(\gamma_j)$ and
$b_j$ a weakly increasing word of length $\gamma_j$ for
$1\le j\le n$.  Fix $1\le i\le n-1$ and
write $b'=\tau_{R_{i+1,R_i}}(b)=b'_n\otimes\dots\otimes b'_1$.
In this case $\tau_{R_{i+1},R_i}$ is an automorphism of
conjugation acting in multiplicity space.
By definition $b'_j=b_j$ for $j\not\in[i,i+1]$
and $P(b'_{i+1} b'_i) = P(b_{i+1} b_i)$.
It follows that $b'$ can be computed from $b$ by
a jeu-de-taquin on the two row skew tableau with word
$b_{i+1} b_i$.
\end{rem}

\subsection{$\slnhat$ crystal structure on rectangular crystals}

Suppose the sequence $R$ consists of a single rectangular partition
with $k$ rows and $l$ columns, so that $B^R=B^{k,l}$.
In \cite{KMN2} the existence of a unique
classical $\slnhat$ crystal structure on $B^{k,l}$ was proved.
The $sl_n$ crystal structure has already been described
in Section \ref{sl_n sec}.  Using the properties of the perfect crystal
$B^{k,l}$ given in \cite{KMN2}, an explicit tableau
construction for $\et_0$ is presented.

The Dynkin diagram of $\slnhat$
admits the rotation automorphism that sends $i$ to $i+1$ modulo $n$.
It follows that there is a bijection
$\psi:B^{k,l}\rightarrow B^{k,l}$ such that
the following diagram commutes for all $i\in I$:
\begin{equation} \label{psi def}
\begin{CD}
  B^{k,l} @>\psi>> B^{k,l} \\
  @V{\ft_i}VV 		@VV{\ft_{i+1}}V \\
  B^{k,l} @>>\psi> B^{k,l}
\end{CD}
\end{equation}
where subscripts are taken modulo $n$.
Of course it is equivalent to require that $\psi$ satisfy the
diagram with $\et_i$ replacing $\ft_i$.

\begin{lem} $\psi$ is unique and 
rotates content in the sense that for
all $b\in B^{k,l}$
\begin{equation} \label{rotate content}
  m_{i+1}(\psi(b))=m_i(b)
\end{equation}
for all $0\le i\le n-1$ where $m_0$ is equal to $m_n$ by convention.
\end{lem}
\begin{proof} Let $b$ be the $sl_n$-highest weight vector in $B^{k,l}$,
given explicitly by $\key((l^k))$.  By the definition
of $\psi$ and the connectedness of $B^{k,l}$ it
is enough to show that $\psi(b)$ is uniquely determined.
By definition $\epsilon_i(b)=0$ for all $i\in J$.
Recall from \cite{KMN2} that for all $b'\in B^{k,l}$,
$\sum_{i\in I} \epsilon_i(b') \ge l$, that $b'$ is said to
be minimal if equality holds, and that for any sequence
of numbers $(a_0,\dots,a_{n-1})\in \N^n$ that sum to $l$,
there is a unique minimal vector $b'$ such that
$\epsilon_i(b')=a_i$ for all $i\in I$.

These facts imply that $\epsilon_0(b)=l$.  
From the definition of $\psi$ it follows that
$\epsilon_i(\psi(b))=l \delta_{i1}.$
Thus $\psi(b)$ is minimal and hence uniquely defined.
\end{proof}

Next it is shown that $\psi$ is uniquely defined by a weaker
condition than \eqref{psi def}.

\begin{lem} \label{psi prop}
$\psi$ is uniquely defined by \eqref{psi def}
for $1\le i\le n-2$ and \eqref{rotate content}.
\end{lem}
\begin{proof} Let $b\in B^{k,l}$, $\bh=b|_{[2,n-1]}$ the restriction
of $b$ to the subinterval $[2,n-1]$, $b'=\psi(b)$,
and $\bh'=b'_{[3,n]}$.  By abuse of notation we shall occasionally
identify a (skew) tableau with its row-reading word.

If $u$ is a word or tableau and $p$ is an integer,
denote by $u+p$ the word or tableau whose entries
are obtained from those of $u$ by adding $p$.
The first goal is to show that $\bh$ and $\bh'-1$ are Knuth-equivalent,
that is, $P(\bh)=P(\bh'-1)$.  By the assumption on $\psi$,
$\bh$ admits a sequence of lowering operators
$e_{i_1} \dots e_{i_p}$ for $2\le i_j\le n-2$, if and only if
$\bh'-1$ admits the sequence $e_{i_1} \dots e_{i_p}$.
Since $\bh$ and $\bh'-1$ are words in the alphabet $[2,n-1]$ this
proves that $P(\bh)=P(\bh'-1)$, by Robinson's
characterization of the RSK map (see section \ref{Kostka sec}).

Now the shape of the tableau $b'$ is a rectangle, so its restriction
$\bh'=b'_{[3,n]}$ to a final subinterval of $[n]$, has antinormal shape.
Hence $\bh'-1$ (and hence $\bh'$) is uniquely determined by $\bh$.  

It only remains to show that the subtableau $b'|_{[1,2]}$ is uniquely
specified.  Its shape must be the partition shape given by the
complement in the rectangle $(l^k)$ with the shape of $b'_{[3,n]}$.
Now $b'|_{[1,2]}$ is a column-strict tableau of partition shape
and contains only ones and twos, so it has at most two rows and
is therefore uniquely determined by its content.  But its content is
specified by \eqref{rotate content}.
\end{proof}

The map $\psi$ is explicitly constructed by exhibiting a map
that satisfies the conditions in Lemma \ref{psi prop}.

The following operation is Sch\"utzenberger's promotion
operator, which was defined on standard tableaux but has
an obvious extension to column-strict tableaux \cite{H} \cite{Schu2}.  
Let $D$ be a skew shape and $b\in\CST(D)$.  The promotion operator
applied to $b$ is computed by the following algorithm.
\begin{enumerate}
\item Remove all the letters $n$ in $b$, which removes from $D$ a horizontal 
strip $H$ (skew shape such that each column contains at most one cell).
\item Slide (using Sch\"utzenberger's jeu-de-taquin \cite{H} \cite{Schu})
the remaining subtableau $b|_{[n-1]}$ to the southeast into the horizontal
strip $H$, entering the cells of $H$ from left to right.
\item Fill in the vacated cells with zeros.
\item Add one to each entry.
\end{enumerate}
The resulting tableau is denoted $\pr(b)\in\CST(D)$ and is called the
promotion of the tableau $b$.

\begin{prop} The map $\psi$ of \eqref{psi def} is given by $\pr$.
\end{prop}
\begin{proof} $\pr$ is content-rotating (satisfies \eqref{rotate content})
and satisfies \eqref{psi def} for $1\le i\le n-2$, since
$\pr(t)|_{[2,n]}=\PSE(t_{[n-1]})+1$ and $\PSE$ commutes with
$sl_n$ crystal operators.  By Lemma \ref{psi prop}, $\pr=\psi$.
\end{proof}

In light of \eqref{psi def}, the operators $\et_0$ and $\ft_0$
on $B^{k,l}$ are given explicitly by
\begin{equation} \label{e_0 def}
\begin{split}
  \et_0 &= \pr^{-1} \circ \et_1 \circ \pr \\
  \ft_0 &= \pr^{-1} \circ \ft_1 \circ \pr.
\end{split}
\end{equation}

\begin{rem} \label{psi rem} Consider again the map $\pr$ on
$b\in B^{k,l}$.  The tableau $\bh=b|_{[n-1]}$ has partition
shape $\la:=(l^{k-1},l-p)$, where $p=m_n(b)$.  Its row-reading word has
Schensted tableau pair $\PPP(\bh)=\bh$ and $\QQQ(\bh)=\key(\la)$.  Let
\begin{equation*}
  \bh':=\pr(b)|_{[2,n]}=\PSE(b|_{[n-1]}),
\end{equation*}
which has antinormal shape and whose complementary shape inside
the rectangle $(l^k)$ must be a single row of length $p$,
that is, $\shape(\bh')=(l^k)/(p)$.  So $\PPP(\bh')=\PPP(\bh)$, and
$\QQQ(\bh')$ is a column-strict tableau of shape $\la$
and content $(l-p,l^{k-1})$, that is, $\QQQ(\bh')=\key(w_0 \la)$
where $w_0$ is the longest element of the symmetric group $S_k$.
So $\PPP(\bh')=\PPP(\bh$ and $\QQQ(\bh')=w_0 \QQQ(\bh)$.
\end{rem}

Let $D$ be the skew shape $R_m\otimes \dots\otimes R_1$,
$B=B^R$ be the tensor product of rectangular crystals,
$b=b_m\otimes \dots\otimes b_1\in B$.
Note that the operator $\pr$ on $\CST(D)=B$ may be described by
\begin{equation*}
  \pr(b)=\pr(b_m)\otimes\dots\otimes \pr(b_1)
\end{equation*}

By the definition of $\et_0$ on a rectangular crystal and the signature
rule, it follows that
\begin{equation} \label{e_0 on rect tensor}
\begin{split}
  \et_0 &= \pr^{-1} \circ \et_1 \circ \pr \\
  \ft_0 &= \pr^{-1} \circ \ft_1 \circ \pr
\end{split}
\end{equation}
as operators on $B=B^R$.

\begin{ex} Let $n=7$, $R=((2,2),(3,3,3),(3,3))$, and $b\in B^R$ given by
the following skew tableau of shape $R_3\otimes R_2\otimes R_1$:
\begin{equation*}
\begin{matrix}
  \times &\times&\times&\times&\times&\times&1&1\\
  \times &\times&\times&\times&\times&\times&2&2\\
  \times&\times&\times&1&1&3& &\\
  \times&\times&\times&2&3&4& &\\
  \times&\times&\times&3&4&5& &\\
  2&4&6& & & & & & \\
  3&5&7& & & & & &
\end{matrix}
\end{equation*}
The element $\pr(b)$ is given by
\begin{equation*}
\begin{matrix}
  \times &\times&\times&\times&\times&\times&2&2\\
  \times &\times&\times&\times&\times&\times&3&3\\
  \times&\times&\times&2&2&4& &\\
  \times&\times&\times&3&4&5& &\\
  \times&\times&\times&4&5&6& &\\
  1&3&5& & & & & & \\
  4&6&7& & & & & &
\end{matrix}
\end{equation*}
The signature for calculating $\et_1$ on $\pr(b)$ is
\begin{equation*}
  \begin{matrix}
    3&2&2&1&1 \\
    -&+&+&+&+
  \end{matrix}
\end{equation*}
So $\et_1$ must be applied to the second tensor position.  Then
$\et_1(\pr(b))$ equals
\begin{equation*}
\begin{matrix}
  \times &\times&\times&\times&\times&\times&2&2\\
  \times &\times&\times&\times&\times&\times&3&3\\
  \times&\times&\times&1&2&4& &\\
  \times&\times&\times&3&4&5& &\\
  \times&\times&\times&4&5&6& &\\
  1&3&5& & & & & & \\
  4&6&7& & & & & &
\end{matrix}
\end{equation*}
Finally $\et_0(b)=\pr^{-1}(\et_1(\pr(b)))$ is given by:
\begin{equation*}
\begin{matrix}
  \times &\times&\times&\times&\times&\times&1&1\\
  \times &\times&\times&\times&\times&\times&2&2\\
  \times&\times&\times&1&3&3& &\\
  \times&\times&\times&2&4&4& &\\
  \times&\times&\times&3&5&7& &\\
  2&4&6& & & & & & \\
  3&5&7& & & & & &
\end{matrix}
\end{equation*}
\end{ex}

\subsection{Action of $\et_0$ on the tableau pair}

In this section an algorithm is given to compute the
tableau pair $(\PPP(\et_0(b)),\QQQ(\et_0(b)))$ of $\et_0(b)$
directly in terms of the tableau pair $(\PPP(b),\QQQ(b))$ of $b$.
In light of \eqref{e_0 on rect tensor} and \eqref{crystal RSK}
with $g=\et_1$, it is enough to give $\PPP(\pr(b))$ and
$\QQQ(\pr(b))$ in terms of $\PPP(b)$ and $\QQQ(b)$.

$\pr(b)=\pr(b_m)\otimes\dots\otimes\pr(b_1)$ can be constructed
by applying Remark \ref{psi rem} to each tensor factor.
The element $b$ is regarded as a skew tableau of shape
$R_m\otimes \dots\otimes R_1$.
Let $\bh=b|_{[n-1]}$ and write
$\bh = \dotsb \otimes \bh_j \otimes \dotsb$, so that
$\bh_j = b_j|_{[n-1]}$.  Let $D_j=(\mu_j^{-1+\eta_j},\mu_j-m_n(b_j))$
denote the shape of $\bh_j$, so that $\bh$ has shape
$D_m \otimes \dots\otimes D_1$.

Next, let $\bh'_j$ be the tableau of skew shape
$D'_j=(\mu_j^{\eta_j})/(m_n(b_j))$ obtained from $\bh_j$
as in Remark \ref{psi rem}.
Write $\bh'=\bh'_m\otimes \dots\otimes \bh'_1$, a skew column-strict tableau
of shape $D'=D'_m\otimes\dots\otimes D'_1$.
Finally, $\pr(b)$ is obtained by adjoining zeros to $\bh'$ at the
vacated positions of the shape $R_m \otimes \dots \otimes R_1$
that are not in $D'$, and then adding one to each entry.

Since $\word(\bh)=\word(b)|_{[n-1]}$, 
\begin{equation*}
  \PPP(\bh)=P(\word(\bh))=P(\word(b)|_{[n-1]})=
  P(\word(b))|_{[n-1]}=\PPP(b)|_{[n-1]}.
\end{equation*}
In other words, $\PPP(\bh)$ is obtained from $\PPP(b)$ by
removing the letters $n$, which occupy a horizontal strip (call it $H$).
It is well-known that $\QQQ(\bh)$ is obtained from $\QQQ(b)$ by
reverse column insertions at the cells of $H$ starting
with the rightmost cell of $H$ and proceeding to the left, ejecting
a weakly increasing word $v$ of length $m_n(b)$.  Another way to say this
is that there is a unique weakly increasing word $v$ of length $m_n(b)$
such that $P(v \QQQ(\bh))=\QQQ(b)$.  So the content of $v$ is the
difference of the contents of $\QQQ(b)$ and $\QQQ(\bh)$.  Since
$\QQQ(b)\in\LRT(R)$, it follows that $v$ is the
weakly increasing word comprised of $m_n(b_j)$ copies of the
maximum letter of $A_j$ for all $j$.

In light of Remarks \ref{coreflection} and \ref{psi rem},
\begin{equation*}
  \QQQ(1+\bh')=\QQQ(\bh') = w_0^R \QQQ(\bh),
\end{equation*}
where $w_0^R$ is the automorphism of conjugation corresponding to the
longest element in the Young subgroup $S_{A_1}\times\dots\times S_{A_m}$
in the symmetric group $S_n$.  Recall that $1+\bh'=\pr(b)|_{[2,n]}$.
Let $H_1$ be the skew shape given by the difference of the shapes of
$\PPP(\pr(b))$ and $\PPP(\pr(b)|_{[2,n]})$.  It is well-known that
$\QQQ(\pr(b)|_{[2,n]})$ is obtained from $\QQQ(\pr(b))$
by reverse row insertions at the cells of $H_1$ starting from the
rightmost and proceeding to the left.  Let $u$ be the ejected word.
Then using an argument similar to that above,
$P(\QQQ(\pr(b)|_{[2,n]}) u)=\QQQ(\pr(b))$ and
$u$ is the weakly increasing word
comprised of $m_n(b_j)$ copies of the minimal letter of $A_j$ for all $j$.
Since $u$ and $v$ are both weakly increasing words it is easy to
calculate directly that $u= w_0^R v$.  Therefore
\begin{equation*}
\begin{split}
  \QQQ(b)&=P(v \QQQ(\bh)) \\
  \QQQ(\pr(b))&=P(\QQQ(1+\bh')u) \\
  &= P((w_0^R \QQQ(\bh))(w_0^R v)).
\end{split}
\end{equation*}

\begin{rem} \label{psi RSK} In summary, the tableau pair
$(\PPP(\pr(b)),\QQQ(\pr(b)))$ is constructed from the tableau pair
$(\PPP(b),\QQQ(b))$ by the following steps.
Let $P=\PPP(b)$ and $Q=\QQQ(b)$.
\begin{enumerate}
\item Let $H$ be the horizontal strip given by the positions of the
letters $n$ in $P$.
\item Let $v$ be the weakly increasing word and
$\Qh$ the tableau such that $\shape(\Qh)=\shape(Q)-H$, such that
$Q=P(v \Qh)$.  These may be produced by reverse column insertions
on $Q$ at $H$ from right to left.
\item Then $\QQQ(\pr(b))=P((w_0^R \Qh)(w_0^R v)$.
\item Let $H_1$ be the horizontal strip $\shape(\QQQ(\pr(b)))-\shape(\Qh)$.
\item Let $P_1$ be the column-strict tableau given by adjoining
to $P|_{[n-1]}$ the letters $n$ at the cells of $H_1$.
\item $\PPP(\pr(b)) = \pr(P_1)$.
\end{enumerate}
By \cite[Proposition 15]{S1},
\begin{equation} \label{chi and e_0}
  (w_0^R \Qh)(w_0^R v)=\chi_R^{-m_n(b)}(v \Qh)
\end{equation}
where $\chi_R$ is the LR analogue of the right circular shift of a
word by positions defined in \cite{S1}.
\end{rem}

All of these steps are invertible, so a description of
$\pr^{-1}$ is obtained as well.

\begin{ex} Continuing the previous example, the
image of $b$ under the map \eqref{R RSK} is given by
the tableau pair $P=\PPP(b)$ and $Q=\QQQ(b)$:
\begin{equation*}
P =
\begin{matrix}
  1&1&1&1\\
  2&2&2&2\\
  3&3&3&3\\
  4&4&4& \\
  5&5& & \\
  6& & & \\
  7& & & 
\end{matrix} \qquad
Q =
\begin{matrix}
  1&1&3&3\\
  2&2&4&6\\
  3&4&5&7\\
  4&5&6& \\
  5&7& & \\
  6& & & \\
  7& & & 
\end{matrix}.
\end{equation*}
Then $H$ is the skew shape given by the single cell $(7,1)$,
$v=7$, $w_0^R v = 6$,
\begin{equation*}
\Qh=
\begin{matrix}
  1&1&3&3\\
  2&2&4&6\\
  3&4&5&7\\
  4&5&6& \\
  5&7& & \\
  6& & &
\end{matrix} \qquad
w_0^R \Qh=
\begin{matrix}
  1&1&3&3\\
  2&2&4&6\\
  3&4&5&7\\
  4&5&6& \\
  5&7& & \\
  7& & &
\end{matrix},
\end{equation*}
and
\begin{equation*}
\QQQ(\pr(b)) =
\begin{matrix}
  1&1&3&3&6\\
  2&2&4&6&\\
  3&4&5&7&\\
  4&5&6& &\\
  5&7& & &\\
  7& & & &
\end{matrix}.
\end{equation*}
So $H_1$ consists of the single cell $(1,5)$ and
\begin{equation*}
P_1 =
\begin{matrix}
  1&1&1&1&7\\
  2&2&2&2\\
  3&3&3&3\\
  4&4&4& \\
  5&5& & \\
  6& & & \\
\end{matrix} \qquad
\PPP(\pr(b)) = \pr(P_1) =
\begin{matrix}
  1&2&2&2&2\\
  3&3&3&3&\\
  4&4&4&4&\\
  5&5&5& &\\
  6&6& & &\\
  7& & & &
\end{matrix}.
\end{equation*}
\end{ex}

\subsection{The $R$-cocyclage and $\et_0$}

In \cite{S1} the $R$-cocyclage relation was defined on
the set $\LRT(R)$.  In the Kostka case this is a weak
version of the dual of the cyclage poset of Lascoux and
Sch\"utzenberger \cite{LS2}.  
It is shown that every covering $R$-cocyclage relation,
realized as recording tableaux, is induced by $\et_0$ on
some element of $B^R$.

\begin{thm} \label{null cocyclage} Let $ux$ be an $R$-LR word
with $x$ a letter.  Then there is an element $b\in B^R$ such that
$\QQQ(b)=P(ux)$ and $\QQQ(\et_0(b))=P(\chi_R(ux))$, provided that
the cell $s=\shape(P(ux))/\shape(P(u))$
is not in the $n$-th row.  In particular, every $R$-cocyclage
covering relation is realized by the action of $\et_0$ in this way.
\end{thm}
\begin{proof} By definition (see \cite{S1}),
every covering relation in the $R$-cocyclage has the form
that $P(v)$ covers $P(\chi_R(v))$ where $v$ is an $R$-LR word.
It follows from \cite[Proposition 23]{S1} that if
$P(ux)\rightarrow P(\chi_R(ux))$ is an $R$-cocyclage then $x>1$.
If $s$ is in the $n$-th row,
then by the column-strictness of $P(ux)$ and the fact that
all letters are in the set $[n]$, $x=1$.  So no $R$-cocyclage
covering relations are excluded by the restriction that
$s$ not be in the $n$-th row.

Let $x'u'=\chi_R(ux)$, that is, $x'=w_0^R x$ and $u'=w_0^R u$.
Since $P(u')=P(w_0^R u)=w_0^R P(u)$ and the automorphism of conjugation
$w_0^R$ preserves shape, without loss of generality it may be
assumed that $u$ is the row-reading word of a column-strict tableau $U$
of partition shape $\lh=\la-s$ where $\la=\shape(\PPP(ux))$.
Now a skew tableau $t$ has partition shape if and only if
$\QQQ(t)=\key(\nu)$, so $u'$ is the word of the column-strict tableau
$U'=w_0^R U$.

Let $Q=P(Ux)$, $\la=\shape(Q)$, $s=(t,\la_t)=\la-\lh$, and
$w=r_1 r_2 \dots r_{t-1}$.  Define $b\in B^R$ by
$\PPP(b)=\key(w\la)$ and $\QQQ(b)=Q$; such a $b$ exists
since $Q\in\LRT(\la;R)$ and \eqref{R RSK} is a bijection.
It must be shown that $\QQQ(\et_0(b))=P(x'U')$.
This shall be verified by applying the formula
\eqref{e_0 on rect tensor} and Remark \ref{psi RSK}.

Now $\la_t>\la_n$ since $s=(t,\la_t)$ is a corner cell and $t<n$.
The horizontal strip $H$ given by the cells of $\key(w\la)$
containing the letter $n$, is merely the $n$-th row of the shape $\la$.
Since $Q\in\CST(\la)$ (and all tableaux are in the alphabet $[n]$),
the subtableau given by the first $\la_n$ columns of $Q$
is equal to $\key((\la_n^n))$.  Let $Q_r$ be the rest of $Q$
and $R'$ the sequence of rectangles obtained by removing the first
$\la_n$ columns from each of the rectangles in $R$.
Since $Q$ is $R$-LR and the column-reading word of $Q$ equals that of
$\key((\la_n^n))Q_r$, it follows that $Q_r$ is $R'$-LR.
Let $y$ be the minimal element of the last interval $A_m$. 
In the notation of Remark \ref{psi RSK}, since $Q_r$ is $R'$-LR, it follows
that $v=n^{\la_n}$, $\Qh=\key((\la_n^{n-1},0))Q_r$, $w_0^R v = y^{\la_n}$,
and $w_0^R \Qh=\key((\la_n^{n-\eta_m},0,\la_n^{-1+\eta_m}))Q_r$.  So
\begin{equation} \label{Q psi}
\begin{split}
  \QQQ(\pr(b)) &= P((w_0^R \Qh)(w_0^R v)) \\
  &= \key((\la_n^{n-\eta_m},0,\la_n^{-1+\eta_m}))Q_r y^{\la_n}.
\end{split}
\end{equation}
The right hand side is a column-strict tableau of shape
$(\la_1+\la_n,\la_2,\dots,\la_{n-1},0)$, so that the horizontal strip
$H_1=(\la_1+\la_n)/(\la_1)$ is entirely in the first row, and  
$\PPP(\pr(b))$ is formed from $1+\key(w(\la_1,\dots,\la_{n-1},0))$
by pushing the first row to the right by $\la_n$ cells and
placing 1's in the vacated positions.  The tableau
$\key(w(\la_1,\dots,\la_{n-1},0))$ contains $\la_t$ ones.
Hence $\PPP(\pr(b))$ contains $\la_t$ twos in its first row,
that is, $\epsilon_1(\PPP(\pr(b)))=\la_t>0$.

So $\PPP(\pr(b))$ admits $\et_0$, which changes the letter
$2$ at the cell $(1,m_n(b)+1)$ to a 1.  By \eqref{crystal RSK},
\begin{equation*}
\begin{split}
  \PPP(\et_1(\pr(b)))&=\et_1(\PPP(\pr(b))) \\
  \QQQ(\et_1(\pr(b)))&=\QQQ(\pr(b)).
\end{split}
\end{equation*}
Now $\pr^{-1}$ is applied, reversing the algorithm in
Remark \ref{psi RSK} starting with the tableau pair
$\PPP(\et_1(\pr(b)))$ and $\QQQ(\et_1(\pr(b)))$,
writing $H_1'$ and $H'$ for the analogues of $H_1$ and $H$.
By Remark \ref{psi RSK} and direct calculation,
\begin{equation} \label{mid P}
\begin{split}
  1 + \PPP(\pr^{-1}(\et_1(\pr(b))))|_{[n-1]} &=
	P(\PPP(\et_1(\pr(b)))|_{[2,n]}) \\
	&=
	1+ \key(w((\la_1,\dots,\la_{n-1},0)-\{s\})).
\end{split}
\end{equation}
In particular $H_1'=H_1\cup\{s\}$.
The reverse row insertions on $\QQQ(\et_1(\pr(b)))=\QQQ(\pr(b))$
at $H_1\subset H_1'$ merely remove the $\la_n$ copies of $y$ from the
first row by \eqref{Q psi}.  The final reverse row insertion (at the cell
$s$) stays within the subtableau $Q_r$ and changes it to the tableau
$Q_r"$ say, and ejects the letter $x$, since 
since $x$ is obtained from $Q$ by reverse row insertion at $s$
and $Q=\key((\la_n^n))Q_r$.  The result of the reverse row insertion
on $\QQQ(\et_1(\pr(b)))$ at $H_1'$ is
\begin{equation*}
  Q_1=\key((\la_n^{n-\eta_m},0,\la_n^{-1+\eta_m}))Q_r",
\end{equation*}
with ejected word $x y^{\la_n}$.  Writing $U'=\key((\la_n^n))Q_r'$
and using the fact that $Q_r$ is $R'$-LR,
\begin{equation*}
\begin{split}
\QQQ(\et_0(b)) &= \QQQ(\pr^{-1}(\et_1(\pr(b)))) \\
&= P((w_0^R x y^{\la_n})
(w_0^R \key((\la_n^{n-\eta_m},0,\la_n^{-1+\eta_m}))Q_r")\\
&= P(x' n^{\la_n} \key((\la_n^{n-1},0)) Q_r')\\
&= P(x' \key((\la_n^n)) Q_r') \\
&= P(x' U').
\end{split}
\end{equation*}
\end{proof}

\begin{rem} Suppose $B^R$ is a Kostka crystal with
$R_j=(\gamma_j)$ for all $1\le j\le n$.  Then
Theorem \ref{null cocyclage} shows that every covering relation
in the cyclage poset $\bigcup_\la \CST(\la;\gamma)$ is realized
in the recording tableaux by $\et_0$ acting on some $b\in B^R$.
\end{rem}

\begin{rem} \label{last column cocyclage} Let $M=\max_i \mu_i$,
the maximum number of columns among the rectangles $R_i$ in $R$.
Suppose $b\in B^R$ is an $sl_n$-highest weight vector
such that $\QQQ(b)$ has shape $\la$, such that $\la_1>M$.
Then $\QQQ(b)\rightarrow \QQQ(\et_0(b))$ is a covering relation
in the $R$-cocyclage.

To see this, let us adopt the notation of the proof of Theorem
\ref{null cocyclage}.  Let $s=(t,\la_t)$ be the corner cell in the last
column of $\la$.  Then $\la_t=\la_1$ so that $w\la=\la$.
By the choice of $b$, $\PPP(b)=\key(\la)=\key(w\la)$.

To apply Theorem \ref{null cocyclage} it must be shown that
$t<n$.  Suppose not.  Then $\la=(\la_1^n)$ and $\QQQ(b)$
is an $R$-LR tableau of shape $(\la_1^n)$.  Since the total
of the heights of the rectangles in $R$ is $n$, it follows that
all of the rectangles in $R$ must have exactly $\la_1$ columns,
contradicting the assumption that $\la_1>M$.

Apply the reverse row insertion on
$\QQQ(b)$ at the cell $s$, obtaining the column-strict tableau
$U$ of shape $\la-\{s\}$ and ejecting the letter $x$.
Then $P(Ux)\rightarrow P(\chi_R(Ux))$ is an $R$-cocyclage,
by \cite[Remark 17]{S1}.
\end{rem}
 
\begin{ex} Continuing the example, $\et_0(b)$ is computed below.
$\QQQ(\et_1(\pr(b)))=\QQQ(\pr(b))$ and
\begin{equation*}
\PPP(\et_1(\pr(b))) =
\begin{matrix}
  1&1&2&2&2\\
  3&3&3&3&\\
  4&4&4&4&\\
  5&5&5& &\\
  6&6& & &\\
  7& & & &
\end{matrix}.
\end{equation*}
Applying $\pr^{-1}$ to the tableau pair of $\et_1(\pr(b))$
and denoting by $P_1'$, $\Qh'$, $v'$ the
corresponding tableaux and word, we obtain
\begin{equation*}
P_1' =
\begin{matrix}
  1&1&1&2&7\\
  2&2&2&3&\\
  3&3&3&7&\\
  4&4&4& &\\
  5&5& & &\\
  6& & & &
\end{matrix} \qquad
w_0^R \Qh' =
\begin{matrix}
  1&1&3&6\\
  2&2&4&7\\
  3&4&5& \\
  4&5&6& \\
  5&7& & \\
  7& & & 
\end{matrix},
\end{equation*}
$w_0^R v' = 36$, $v'=57$,
\begin{equation*}
\Qh' =
\begin{matrix}
  1&1&3&6\\
  2&2&4&7\\
  3&3&5& \\
  4&4&6& \\
  5&7& & \\
  6& & & \\
   & & & 
\end{matrix}\qquad
\QQQ(\et_0(b)) =
\begin{matrix}
  1&1&3&6\\
  2&2&4&7\\
  3&3&5& \\
  4&4&6& \\
  5&5&7& \\
  6& & & \\
  7& & & 
\end{matrix},
\end{equation*}
and finally
\begin{equation*}
\PPP(\et_0(b)) =
\begin{matrix}
  1&1&1&2\\
  2&2&2&3\\
  3&3&3& \\
  4&4&4& \\
  5&5&7& \\
  6& & & \\
  7& & & 
\end{matrix} \qquad
\end{equation*}
\end{ex}

\begin{rem} Suppose that $R=(R_1,\dots,R_m)$ is such that
each $R_j$ has a common number of columns, say $l$.
Then $B^R$, being the tensor product of perfect crystals
$B^{\eta_j,l}$ of level $l$, is perfect of level $l$
and therefore connected by \cite{KMN}.  In this case, more is true.
Using Remark \ref{last column cocyclage}, every element can be connected
to the unique $sl_n$-highest weight vector in $B^R$ of charge zero,
by applying last column $R$-cocyclages on the recording tableau.

For $R$ a general sequence of rectangles, $B^R$ is still connected.
However it is not necessarily possible to use $\et_0$ to connect
every $sl_n$-component to the zero energy component in such a way
that the energy always drops by one.  For example, take
$n=3$, $R=((2),(1),(1)$ and the $sl_3$-component with $\QQQ$-tableau
$\begin{matrix}1&1\\2&3\end{matrix}$.  The applications of $\et_0$ on
the five elements of this component that admit $\et_0$,
all produce elements in the component with
$\QQQ$-tableau $\begin{matrix}1&1&3\\2& & \end{matrix}$ which
has the same energy.
\end{rem}

\subsection{Rectangle-switching bijections and $\et_0$}

\begin{prop} Let $R_1$ and $R_2$ be rectangles.  Then
the rectangle-switching bijection
\begin{equation*}
\sigma_{R_2,R_1}: B^{R_2}\otimes B^{R_1} \rightarrow
  B^{R_1}\otimes B^{R_2}
\end{equation*}
is an isomorphism of classical $\slnhat$ crystals.
\end{prop}
\begin{proof} Since $\sigma=\sigma_{R_2,R_1}$ is known to be an isomorphism
of $sl_n$-crystals, it only remains to show that $\sigma$ commutes
with $\et_0$.  Let $b\in B^{R_2}\otimes B^{R_1}$.  By the bijectivity
of \eqref{R RSK} it is enough to show that
\begin{equation} \label{P and Q}
\begin{split}
\PPP(\et_0(\sigma(b)))&=\PPP(\sigma(\et_0(b)))\\
\QQQ(\et_0(\sigma(b)))&=\QQQ(\sigma(\et_0(b))).
\end{split}
\end{equation}
Consider first the process in passing from $(\PPP(b),\QQQ(b))$ to
$(\PPP(\et_0(b)),\QQQ(\et_0(b)))$.  Let $H$ and $H_1$ be the
be the horizontal strips, 
$v$ the weakly increasing word and $\Qh$ the tableau as in
Remark \ref{psi RSK}.  Let $b'=\sigma(b)$, and let
$H'$, $H_1'$, $v'$, $\Qh'$ the analogous objects in passing from
$(\PPP(b'),\QQQ(b'))$ to $(\PPP(\et_0(b')),\QQQ(\et_0(b')))$.

Let $\la=\shape(\PPP(b))$.
Observe that $\PPP(b')=\PPP(\sigma(b))=\PPP(b)$
by \eqref{rect switch}, so that $\shape(\PPP(b'))=\la$
and $H'=H$.  This implies that
the increasing words $v'$ and $v$ have the same length
$m_n(b)$ and $\shape(\Qh')=\shape(\Qh)$; call their common shape
$\lh$.  Now $P(v' \Qh')=\QQQ(b')$ is the unique element in
$\LRT(\la;(R_2,R_1))$ and $P(v \Qh)=\QQQ(b)$ is the unique element in
$\LRT(\la;(R_1,R_2))$.  So
\begin{equation} \label{P tau prime}
\begin{split}
  P(\tau(v \Qh))&=\tau P(v \Qh))\\
  &= \tau \QQQ(b) = \QQQ(b') = P(v' \Qh').
\end{split}
\end{equation}
On the other hand, consider the word $v\Qh$, identifying $\Qh$ with its
row-reading word.  The tableau $Q(v \Qh)$ has shape $\la$.
Let $p$ and $q$ be the number of cells in $\Qh$ and $Q$ respectively.
Since $\Qh$ is a column-strict tableau of shape $\lh$,
$Q(v\Qh)|_{[p]}$ is the rowwise standard tableau of shape $\lh$,
the unique standard tableau of shape $\lh$ in which $i+1$ is located
immediately to the right of $i$ provided $i+1$ is not in the first column.
Now $Q(v\Qh)|_H$ is filled from left to right by the numbers $p+1$ through
$q$, since it records the column-insertion of the weakly increasing word $v$.
The same argument applies to $Q(v'\Qh')$, so that $Q(v\Qh)=Q(v'\Qh')$.
By this, \eqref{P tau prime}, and \eqref{tau-prime}, it follows that
$\tau(v \Qh)=v' \Qh'$.  Write $R'=(R_2,R_1)$, $w_0^{R'}$ and
$\chi_{R'}$ for the corresponding constructions for $R'$.  By
\cite[Theorem 16]{S1},
\begin{equation*}
\begin{split}
  (w_0^{R'} \Qh')(w_0^{R'} v') &=\chi_{R'}^{-m_n(b)}(v' \Qh')\\
  &= \chi_{R'}^{-m_n(b)}(\tau(v \Qh)) \\
  &= \tau(\chi_R^{-m_n(b)}(v \Qh)) \\
  &= \tau((w_0^R \Qh)(w_0^R v)).
\end{split}
\end{equation*}
Applying the $P$ tableau part of \eqref{tau-prime} to the word
$(w_0^R \Qh)(w_0^R v)$,
\begin{equation*}
\begin{split}
\tau(\QQQ(\et_0(b)))&=\tau(P((w_0^R \Qh)(w_0^R v)))\\
  &= P(\tau((w_0^R \Qh)(w_0^R v)))\\
  &= P((w_0^{R'} \Qh')(w_0^{R'} v'))\\
  &= \QQQ(\et_0(b')).
\end{split}
\end{equation*}
By this and the $\QQQ$ tableau part of \eqref{rect switch}
for the word $\et_0(b)$,
\begin{equation*}
\begin{split}
\QQQ(\sigma(\et_0(b)))&=\tau(\QQQ(\et_0(b)))\\
&=\QQQ(\et_0(b'))\\
&=\QQQ(\et_0(\sigma(b))
\end{split}
\end{equation*}
This proves the equality of the $\QQQ$-tableaux in \eqref{P and Q}.

For the $\PPP$-tableaux, let us recall the process that leads from
$(\PPP(b),\QQQ(b))$ to
$(\PPP(\et_0(b)),\QQQ(\et_0(b)))$ and from
$(\PPP(b'),\QQQ(b'))$ to
$(\PPP(\et_0(b')),\QQQ(\et_0(b')))$ where $b'=\sigma(b)$.
Recall that $\PPP(b')=\PPP(b)$.  Clearly they have the same
restriction to the alphabet $[n-1]$.
On the other hand, since it has been shown that the tableaux
$\QQQ(\et_0(b))$ and $\QQQ(\et_0(b'))$ have the same shape,
it follows that the horizontal strips $H_1$ and $H_1'$
coincide.  So $\PPP(\et_0(b))$ and $\PPP(\et_0(b'))$ coincide when
restricted to the alphabet $[2,n]$.  Since both tableaux also have the same
partition shape they must coincide.  This, together
with the $\PPP$-tableau part of \eqref{rect switch} applied to $\et_0(b)$,
shows that
\begin{equation*}
\begin{split}
  \PPP(\et_0(\sigma(b)))&=\PPP(\et_0(b'))\\
  &= \PPP(\et_0(b))=\PPP(\sigma(\et_0(b))).
\end{split}
\end{equation*}
\end{proof}

\subsection{Energy function}

In this section it is shown that the energy function of the
classical $\slnhat$ crystal $B^R$ is given by the
generalized charge of \cite{S1} on the $\QQQ$-tableau.
The definition of energy function follows \cite{KMN2}
and \cite{NY}.

Consider the unique classical $\slnhat$ crystal isomorphism
\begin{equation*}
\begin{split}
  \sigma: B^{R_2} \otimes B^{R_1} &\cong B^{R_1} \otimes B^{R_2} \\
  b_2 \otimes b_1 &\mapsto b_1' \otimes b_2'.
\end{split}
\end{equation*}

An energy function $H:B^{R_2}\otimes B^{R_1}\rightarrow \Z$
is a function that satisfies the following axioms:
\begin{enumerate}
\item[(H1)] $H(\ft_i(b))=H(b)$ for all $i\in J$ and
$b\in B^{R_2}\otimes B^{R_1}$
such that $\ft_i(b)$ is defined, and similarly for $\et_i$.
\item[(H2)] For all $b=b_2\otimes b_1\in B^{R_2}\otimes B^{R_1}$
such that $\et_0(b)$ is defined,
\begin{equation*}
  H(\et_0(b)) - H(b) =
  \begin{cases}
     1 & \text{if $\epsilon_0(b_2)\le \phi_0(b_1)$ and
	$\epsilon_0(b'_1)\le\phi_0(b'_2)$} \\
    -1 & \text{if $\epsilon_0(b_2)>\phi_0(b_1)$ and
	$\epsilon_0(b'_1)>\phi_0(b'_2)$} \\
     0 & \text{otherwise.}
  \end{cases}
\end{equation*}
\end{enumerate}

If $H':B^{R_1}\otimes B^{R_2}\rightarrow \Z$ is defined in the same way
then $H'\circ \sigma = H$.

\begin{lem} \label{connected} $B^{R_2}\otimes B^{R_1}$ is connected.
\end{lem}
\begin{proof} Without loss of generality assume that
$\mu_1\ge\mu_2$.  Let
\begin{equation*}
\gamma(R)=(\mu_1^{\eta_1},\mu_2^{\eta_2})
\end{equation*}
be the partition whose Ferrers diagram is obtained by placing the
shape $R_1$ atop $R_2$.  Let $Q_R$ be the unique element of the
singleton set $\LRT(\gamma(R);(R_1,R_2))$.  Observe that
$\gamma(R)$ is the only shape admitting an $R$-LR tableau
that has at most $\mu_1$ columns.
Let $v_R\in B^R$ be the unique $sl_n$-highest weight vector of weight
$\gamma(R)$; it is given explicitly by $v_R = Y_2 \otimes Y_1$
in the notation of the definition of $R$-LR, and satisfies
$\PPP(v_R)=\key(\gamma(R))$ and $\QQQ(v_R)=Q_R$.  

It is shown that every element $b\in B^R$ is connected to $v_R$.
First, using $sl_n$-raising operators, it may be assumed that
$b$ is an $sl_n$-highest weight vector.  If $\QQQ(b)$ has at
most $\mu_1$ columns then $b=v_R$.  Otherwise $\QQQ(b)$ has more
than $\mu_1$ columns, and Remark \ref{last column cocyclage} applies.
But $\et_0(b)$ is closer to $v_R$ in the sense that
$\QQQ(\et_0(b))$ has one fewer cells to the right of the $\mu_1$-th
column than $\QQQ(b)$ does \cite[Proposition 38]{S1}, so induction
finishes the proof.
\end{proof}

By Lemma \ref{connected} $H$ is uniquely determined up to a global
additive constant.  $H$ is normalized by the condition
\begin{equation*}
  H(v_R) = 0
\end{equation*}
with $v_R$ as in Lemma \ref{connected}.  Equivalently,
\begin{equation*}
  H(\key(R_2) \otimes \key(R_1)) = \min(\eta_1,\eta_2) \min(\mu_1,\mu_2),
\end{equation*}
where $\key(R_i)$ is the highest-weight vector in $B^{R_i}$ and
the value of $H$ is the size of the rectangle $R_1\cap R_2$.

For a tableau $Q\in\LRT(\la;(R_1,R_2))$, define
$d_{R_1,R_2}(Q)$ to be the number of cells in the shape of
$Q$ that are strictly east of the $\max(\mu_1,\mu_2)$-th column.

\begin{prop} \label{two energy}
Let $R=(R_1,R_2)$ be a pair of rectangles.  Then for all
$b\in B^R$, $H(b)=d_{R_1,R_2}(\QQQ(b))$.
\end{prop}
\begin{proof} Follows immediately from the proof of Lemma \ref{connected}.
\end{proof}

Now consider $R=(R_1,\dots,R_m)$ with $B_j=B^{R_j}$ and
$B=B^R=B_m\otimes\dots\otimes B_1$.
The energy function for $B$ is given as follows.  Denote by
$H_{i,w}: w B \rightarrow \Z$ given by the value of the
energy function $H_{(wB)_{i+1},(wB)_i}$ at the $(i+1)$-st and $i$-th
tensor positions (according to our convention).
Recall the isomorphisms of classical $\slnhat$-crystals \eqref{sigma perm}.
Define the cyclic permutation $w_{i,j}=r_{i+1} r_{i+2} \dotsb r_{j-1}$
for $1\le i<j\le m$.  Define $E_R: B\rightarrow \Z$ by
\begin{equation*}
  E_R(b) = \sum_{1< j\le m} \sum_{1\le i< j}
	H_{i,w_{i,j}}(\sigma_{R,w_{i,j}R}(b)).
\end{equation*}
Call the inner sum $E_R^{(j)}(b)$.  

The following version of \cite[Lemma 5.1]{KMOTU} holds for $E_R$
with no additional difficulty.

\begin{lem} \label{energy and e} Let $B=B^R$ where
$R_j=(\mu_j^{\eta_j})$ and $l\ge\mu_j$ for all $j$.
Suppose $\epsilon_0(b)>\mu_j$ for all $j$.  Write
$\et_0(b) = \dotsb \otimes b_{k+1} \otimes \et_0(b_k) \otimes b_{k-1} \dotsb$.
Then $k>1$ and $E_R^{(j)}(\et_0(b))=E_R^{(j)}(b)-\delta_{j,k}$.
\end{lem}

\begin{ex} Let $b$ be as in the running example.  Then
\begin{equation*}
  \tau_2(b)=
\begin{matrix}
  \times &\times&\times&\times&\times&\times&1&1\\
  \times &\times&\times&\times&\times&\times&2&2\\
  \times&\times&\times&1&3&3& & \\
  \times&\times&\times&2&4&4& & \\
  1&3&5& & & & & \\
  2&4&6& & & & & \\
  3&5&7& & & & &
\end{matrix}
\end{equation*}
Write $\tau_2(b)=b_2'\otimes b_3'\otimes b_1$.  Now
$P(\word(b_2\otimes b_1))$ has shape $(4,3,3,2,1)$ so that
$d_1(b)=1$.  $P(\word(b_3 \otimes b_2))$ has shape $(4,4,3,2,2)$ so that
$d_2(b)=2$.  Finally $P(\word(b_3'\otimes b_1))$ has shape $(3,3,2,2)$
so that $d_1(\tau_2(b))=0$.  So $E_R(b)=1+2+0=3$.
\end{ex}

\subsection{Energy and generalized charge}

Define the map $E_R:\LRT(R)\rightarrow\Z$ by
$E_R(Q)=E_R(b)$ for any $b\in B^R$ such that $\QQQ(b)=Q$.
This map is well-defined since $E_R$ is constant
on $sl_n$-components and the map \eqref{R RSK} is a bijection.
It follows immediately from the definitions that
\begin{equation*}
  E_R(Q) = \sum_{1\le i<j\le m} d_{R_i,R_j}(\tau_{R_{i+1},R_j}\circ
  \dots \circ \tau_{R_{j-1},R_j}(Q).
\end{equation*}

The Kostka case of the following result was first
proven by K. Kilpatrick and D. White.  In the further special
case that $\mu$ is a partition it was shown in \cite{NY}
that $E_R(Q)$ is the charge.  Now in the Kostka case the generalized charge
statistic $\charge_R$ specializes to the formula of charge
in \cite{LLT}.

\begin{thm} \label{charge is energy} $\charge_R=E_R$ on $\LRT(R)$.
\end{thm}
\begin{proof} Let $Q\in\LRT(R)$ of shape $\la$, say.  It
will be shown by induction on the number of rectangles $m$ and then on
$\charge_R(Q)$, that
that $E_R:\LRT(R)\rightarrow\Z$ satisfies the intrinsic characterization
of $\charge_R$ by the properties (C1) through (C4) \cite[Theorem 21]{S1}.
Let $M=\max_j \mu_j$.

First, (C2) need only be checked when $Q'<_R Q$ is a last column
$R$-cocyclage, and (C4) need only be verified when
$\la_1=M$.  To see this, if $\la_1>M$ then
there is a last column $R$-cocyclage $Q'<_R Q$ and in this case
$\charge_R(Q')=\charge_R(Q)-1$.  Otherwise $\la_1=M$.
If (C3) does not apply, then one may apply (C4) several times to
switch the widest rectangle closer to the beginning of
the sequence $R$ and then apply (C3), which decreases the number of
rectangles $m$.

(C1) is trivial.  For (C2),
let $Q' <_R Q$ be a last-column $R$-cocyclage with $\shape(Q)=\la$
such that $\la_1>M$.  Let $P=\key(\la)$ and
$b$ such that $\PPP(b)=\key(\la)$ and $\QQQ(b)=Q$ as in
Theorem \ref{null cocyclage} and
Remark \ref{last column cocyclage}.  By the proof of 
Theorem \ref{null cocyclage} in this case,
$\epsilon_0(b)=\epsilon_1(\PPP(\pr(b)))=\la_t=\la_1>M$.
By Lemma \ref{energy and e},
\begin{equation*}
E_R(\QQQ(\et_0(b)))=E_R(\et_0(b))=E_R(b)-1=E_R(\QQQ(b))-1=E_R(Q)-1.
\end{equation*}
However, $\QQQ(\et_0(b))=Q'$ by Remark \ref{last column cocyclage},
so $E_R(Q')=E_R(Q)-1$, and (C2) has been verified.

To check (C3), let $\Rh=(R_2,\dots,R_m)$ and $\Qh=Q|_{[\eta_1+1,n]}$.
Then $Q$ consists of $\key(R_1)$ sitting atop $\Qh$
and $\Qh\in\LRT(\Rh)$.  It follows that
\begin{equation*}
  d_{R_1,R_j}(\tau_{R_2,R_j} \tau_{R_3,R_j}\dots \tau_{R_{j-1},R_j} Q) = 0
\end{equation*}
for all $j>1$.  Therefore
\begin{equation*}
\begin{split}  
 E_R(Q) &= \sum_{1\le i<j\le m}
   d_{R_i,R_j}(\tau_{R_{i+1},R_j} \dots \tau_{R_{j-1},R_j} Q) \\
&=\sum_{2\le i<j\le m}
   d_{R_i,R_j}(\tau_{R_{i+1},R_j} \dots \tau_{R_{j-1},R_j} Q) \\
&= E_{\Rh}(\Qh), 
\end{split}  
\end{equation*}
which verifies (C3).

For (C4), the proof may be reduce to the case $m=3$.
By abuse of notation we suppress the notation for the sequence
of rectangles, writing $\tau_p(Q)$ for the operator
that acts on the restriction of an LR tableau to the
$p$-th and $(p+1)$-st alphabets, and similarly for the function $d_p$.
Write $w_{i,j} := \tau_{i+1} \tau_{i+2} \dots \tau_{j-1}$
for $1\le i<j\le m$.

Fix $1\le p \le m-1$ and $1\le i<j\le m$.  Write
\begin{equation*}
\begin{split}
  d'_{i,j} &:= d_i(w_{i,j} \tau_p Q) \\
  d_{i,j} &:= d_i(w_{i,j} Q).
\end{split}
\end{equation*}
The value $d'_{i,j}$ is computed using a case by case analysis.
\begin{enumerate}
\item $i<p+1$.  In this case it is clear that $d_{i,j}'=d_{i,j}$.
\item $i=p+1$.  Then $w_{p+1,j}\tau_p=\tau_p w_{p+1,j}$ and 
\begin{equation*}
  d_{p+1,j}' = d_{p+1}(\tau_p w_{p+1,j} Q).
\end{equation*}
\item $p=i$ and $p+1<j$.  Then
\begin{equation*}
  w_{p,j}\tau_p = \tau_{p+1} w_{p+1,j} \tau_p = \tau_{p+1} \tau_p w_{p+1,j}
\end{equation*}
so that
\begin{equation*}
  d_{p,j}' = d_p(\tau_{p+1}\tau_p w_{p+1,j} Q).
\end{equation*}
\item $p=i$ and $p+1=j$.  Here $w_{i,j}$ is the identity, and
\begin{equation*}
  d_{p,p+1}' = d_p(\tau_p Q)=d_p(Q)=d_{p,p+1}.
\end{equation*}
\item $i<p$ and $p+1<j$.  Then $w_{i,j} \tau_p = \tau_{p+1} w_{i,j}$ so that
\begin{equation*}
  d_{i,j}' = d_i(\tau_{p+1} w_{i,j} Q) = d_i(w_{i,j} Q)=d_{i,j}.
\end{equation*}
since the restriction of $w_{i,j}Q$ to the $i$-th and $i+1$-st subalphabets
is not affected by $\tau_{p+1}$.
\item $i<p$ and $j=p+1$.  Then $w_{i,p+1}\tau_p = w_{i,p}$ and
\begin{equation*}
  d_{i,p+1}' = d_i(w_{i,p} Q) = d_{i,p}.
\end{equation*}
\item $i<p$ and $j=p$.  Then $w_{i,p}\tau_p = w_{i,p+1}$ and
\begin{equation*}
  d_{i,p}' = d_i(w_{i,p+1} Q)= d_{i,p+1}.
\end{equation*}
\item $j<p$.  In this case it is clear that $d_{i,j}=d'_{i,j}$.
\end{enumerate}
Based on these computations, the difference in energies
$E_{\tau_p R}(\tau_p Q)-E_R(Q)$ is given as follows.
In cases 1, 4, and 5, and 8, $d'_{i,j}=d_{i,j}$ so these terms
cancel.  The sum of the terms in cases 6 and 7 cancel.  So it is
enough to show that the sum of terms in 2 and 3 cancel, that is,
\begin{equation*}
\begin{split}
 0 &= \sum_{j>p+1} (d_{p+1}(\tau_p w_{p+1,j} Q)-d_{p+1}(w_{p+1,j}Q)) \\
&+ \sum_{j>p+1} (d_p(\tau_{p+1}\tau_p w_{p+1,j} Q)-d_p(w_{p,j}Q)).
\end{split}
\end{equation*}
Rewriting $d_p(w_{p,j}Q)=d_p(\tau_p w_{p+1,j}Q)$, observe that
without loss of generality it may be assumed that $m=3$ and it must
be shown that
\begin{equation} \label{three rectangles}
  0 = d_2(\tau_1\tau_2 Q)-d_2(Q)+d_1(\tau_2 \tau_1 Q)-d_1(Q).
\end{equation}

Recall that in verifying (C4) it may be assumed 
that $\la_1=M$.  In this case \cite[Remark 39]{S1} applies.
Say $\mu_k=\la_1$.  Then $d_k(Q)=0$ and if $k>1$, $d_{k-1}Q=0$.
There are three cases, namely $k=1$, $k=2$, or $k=3$.  If $k=3$ then
all four terms in \eqref{three rectangles} are zero.
If $k=1$ then the first and fourth terms are zero and the
second and third agree, while if $k=2$ then the
second and third are zero and the first and fourth agree.
\end{proof}

\begin{cor} \label{LRT and path}
\begin{equation*}
  \sum_{b\in B^R} e^{\wt(b)} q^{E_R(b)} =
  \sum_\la \ch V^{\slwt(\la)} K_{\la;R}(q),
\end{equation*}
where $\la$ runs over partitions of length at most $n$.
\end{cor}
\begin{proof} The equality follows immediately from
Theorem \ref{charge is energy},
the weight-preserving bijection \eqref{R RSK}, and
\cite[Theorem 11]{S1}.
\end{proof}

\section{Tensor product structure on Demazure crystals}

The tensor product structure for the Demazure crystals,
is a consequence of an inhomogeneous version of
\cite[Theorem 2.3]{KMOU} that uses Lemma \ref{energy and e}.

\begin{thm} \label{crystal tensor} Let $\mu=(\mu_1,\dots,\mu_m)$ be
a partition of $n$.  Then
\begin{equation*}
  \BB_{w_\mu}(l \La_0) \cong 
  B^{\mu_m,l}\otimes  \dots\otimes B^{\mu_1,l} \otimes u_{l\La_0}
\end{equation*}
as classical $\slnhat$ crystals, where the affine $\slnhat$ Demazure crystal
is viewed as a classical $\slnhat$-crystal by composing its weight
function $\wt$ with the projection $P\rightarrow\PCL$.  Moreover,
if $v\mapsto b\otimes u_{l\La_0}$ then 
\begin{equation*}
  \inner{d}{l\La_0-\wt(v)} = E_R(b)
\end{equation*}
where the left hand side is the distance along the null root $\delta$
of $v$ from the highest weight vector $u_{l\La_0}\in\VV(l\La_0)$ and
$R$ is defined by $R_j=(l^{\mu_j})$ for $1\le j\le m$.
\end{thm}

\section{Proof of Theorem \ref{main result}}

Theorem \ref{main result} follows from Theorem \ref{crystal tensor}
and Corollary \ref{LRT and path}.

\section{Generalization of Han's monotonicity for Kostka-Foulkes
polynomials}

The following monotonicity property for the Kostka-Foulkes
polynomials was proved by G.-N. Han \cite{Ha}:
\begin{equation*}
  \K_{\la,\mu}(q) \le \K_{\la\cup\{a\},\mu\cup\{a\}}(q) 
\end{equation*}
where $\la\cup\{a\}$ denotes the partition obtained by
adding a row of length $a$ to $\la$.

Here is the generalization of this result
for the polynomials $\K_{\la;R}(q)$ that was conjectured
by A. N. Kirillov.

\begin{thm} Let $R$ be a dominant sequence of rectangles and $(k^m)$
another rectangle.  Then
\begin{equation*}
  \K_{\la;R}(q) \le \K_{\la\cup (k^m);R\cup (k^m)}(q)
\end{equation*}
where $\la\cup (k^m)$ is the partition obtained by adding $m$ rows
of size $k$ to $\la$ and $R \cup (k^m)$ is any dominant sequence
of rectangles obtained by adding the rectangle $(k^m)$ to $R$.
\end{thm}
\begin{proof} Write $R=(R_1,\dots,R_t)$, $R_0=(k^m)$, and
$R^+=(R_0,R_1,\dots,R_t)$. 
Define the map
$i_R:\LRT(\la;R)\rightarrow \LRT(\la\cup(k^m);R^+)$ by
$i_R(Q) = P((Q+m) Y_0)$ where $Y_0=\key((k^m))$.
Since the letters of $Y_0$ are smaller than those of $Q+m$, it follows
that $\shape(i_R(Q))=\la\cup(k^m)$.  Moreover $i_R(Q)$ is $R^+$-LR since
it is Knuth equivalent to a shuffle of $Y_0$ and the tableau
$Q+m$, which is $R$-LR in the alphabet $[m+1,n+m]$.
Thus the map $i_R$ is well-defined.  Let $B$ represent the union of the
zero-th and first subalphabets for $w_{0,j} R^+$ and let $Y$ be
the key tableau for the first subalphabet of $w_{0,j}R^+$.  Then
\begin{equation} \label{zero and j}
\begin{split}
  d_0(w_{0,j} i_R(Q)) &= d_0(w_{0,j} Q Y_0) = d_0((w_{0,j}Q) Y_0) \\
  &= d_0(((w_{0,j}Q) Y_0)|_B = d_0(Y Y_0) = 0,
\end{split}
\end{equation}
by the Knuth invariance of $d_0$, the fact that $w_{0,j}$ doesn't
touch letters in the zero-th subalphabet, the definition of $d_0$,
the fact that $w_{0,j} Q\in \LRT(w_{0,j}R)$, and direct calculation
of the shape of $P(Y Y_0)$ combined with Proposition \ref{two energy}.
If $i>0$ then
\begin{equation} \label{i and j}
\begin{split}
  d_i(w_{i,j} i_R(Q)) &= d_i((w_{i,j} i_R(Q))|_{[m+1,m+n]} \\
  &= d_i(w_{i,j} (i_R(Q)|_{[m+1,m+n]})) \\
  &= d_i(w_{i,j} Q).
\end{split}
\end{equation}
From \eqref{zero and j} and \eqref{i and j} it follows that
$E_{R^+}(i_R(Q)) = E_R(Q)$.
\end{proof}

\end{document}